## Sturm-Liouville Theory and Orthogonal Functions

H. Azad and M. T. Mustafa
Department of Mathematics and Statistics,
King Fahd University of Petroleum & Minerals, Dhahran, Saudi Arabia
hassanaz@kfupm.edu.sa, tmustafa@kfupm.edu.sa

ABSTRACT. We revisit basics of classical Sturm-Liouville theory and, as an application, recover Bochner's classification of second order ODEs with polynomial coefficients and polynomial solutions by a new argument. We also outline how a wider class of equations with polynomial solutions can be obtained by allowing the weight to become infinite at isolated points.

For higher order equations, we also give the basic analysis required for determining the weight functions and constraints on the coefficients which make the differential operator defined by the equation self adjoint for even orders and anti self adjoint in odd orders. We also give explicit examples of such equations.

## 1. INTRODUCTION

In this paper we revisit the basics of Sturm-Liouville theory- for all orders. As an application we determine all differential equations of order 2 which have polynomial coefficients and polynomial solutions.

Problems of this type were first considered by Bochner [3]. A review of this paper is available in the Zentralblatt database [11]. The approach taken in this paper is quite different, and this is explained in detail in Section 3, where the reader will also find a summary of Bochner's method, extracted from the review in Zentralblatt.

It is the duty of every generation to make mathematics as transparent as it can for its efficient transmission to the next generation. It is in this spirit that we have seen fit to put this paper in the public domain.

The Sturm-Liouville systems are equations of the type $(pv')' - qv + \lambda\rho(v) = 0$, where the functions $p$ and $\rho$ are, differentiable and positive in an interval. The solution $v$ is required to satisfy boundary conditions of the type $\alpha_1 v + \beta_1 v' = 0$, and $\alpha_2 v + \beta_2 v' = 0$ at the end points of the interval where $\alpha_1, \beta_1$ are not both zero and $\alpha_2, \beta_2$ are not both zero.

There are variations of these boundary conditions for infinite intervals, or where the functions $p$ and $\rho$ are not positive [8, p. 291], [9, p. 107-109].

Any equation of the type $a(x)y'' + b(x)y' + c(x)y = \lambda y$ can always be written in the above form by multiplying by a suitable weight function [2, p. 45].

A natural question is: what is the explanation for the weight function and the particular form of the boundary conditions.

In this paper, which is also meant to be a contribution to the teaching of these ideas, we show that both the weight and general boundary conditions are forced upon us as soon as we demand that the operator $L(y) = a(x)y'' + b(x)y' + c(x)y$ should be self-adjoint for some weight function $p$.

As an application of this, we recover Bochner's classification of second order operators $L$ which map polynomials of order $n$ to polynomials of order $n$ - for which the natural weight function associated to $L$ is positive- except for isolated zeroes- in the interval (finite or infinite) and for which polynomials have finite norms.

We also outline how one can recover a wider class of equations which have polynomial solutions by allowing the weight to become infinite at isolated points. The Jacobi differential equation is of this type- as explained in sections 2 and 3.

Since we want $L$ to operate on all polynomials of order $n$ there should be no boundary conditions. From the formulas in Section 2, no boundary conditions are needed at a point where the leading term $a(x)$ is zero or where the weight function vanishes.

For example, if the interval is finite, and the leading term vanishes at, say, two points, we may assume that that this interval is [-1, 1]. The finiteness of the weight function at the end points determines the differential equation. So, for each degree n, there will be, up to a constant, a unique polynomial which will be an eigen function of L, and the corresponding eigenvalues would therefore be determined from the form of the equation. Since these eigenvalues are distinct for different degrees, these polynomials would automatically be orthogonal: see Section 3 for details.

Here is an outline of this paper. In Section 2, we revisit the formal aspects of Sturm-Liouville theory and see how the classical orthogonal polynomials fit in this framework. Section 3 deals with canonical forms of $2^{nd}$ order equations whose eigenfunctions are polynomials of finite norm.

In the final section, we write down the determining equations of Sturm-Liouville type of low orders: there are similar formulas for all orders [1]. It turns out that there are no self-adjoint odd order equations- only anti- self adjoint equations. However, if one looks for self adjoint operators in complex differential equations, then one has such equations of all orders. The key to understanding all of this is the case of order 1 equations. We discuss this case briefly, leaving the details for this and higher order equations to a paper under preparation.

After writing the first draft of this paper, we found the paper [7], where we also found the basic reference [3] .The subject of polynomial solutions seems to be undergoing a revival of sorts, as is clear from the bibliography of [7].

From the bibliography of [6], the classification of differential equations which have polynomial solutions is also discussed in several papers, notably by Feldmann [6], Lesky [10] and Everitt et al [4,5].

## 2. BASICS OF STURM-LIOUVILLE THEORY (REVISITED)
## AND CLASSICAL EXAMPLES

In this section, we first give a proof of the following well-known result (cf. [8, p. 277-280] and [2, p. 42-45]). To avoid trivial refinements, the functions appearing here will be assumed to be infinitely differentiable in suitable intervals.

**Proposition** : *Let $L$ be the operator defined by $L(y) = a(x)y'' + b(x)y' + c(x)y$*

*If $p(x) = e^{\int (\frac{b(x)-a'(x)}{a(x)})dx}$ is finite over an interval $I$ and $C$ is a linear space of functions for which the following hold:*

*(i) $C$ is invariant under $L$*

*(ii) for every function $y$ in $C$ , the integral of $py^2$ over $I$ is finite*

*(iii) for any functions $u$ and $y$ in $C$ , the difference of*

$$p(x)a(x)(u'(x)y(x) - u(x)y'(x)) = pa(u'y - uy'),$$

*evaluated at the end-points of $I$ vanishes , then the operator $L$ is self adjoint on $L$ with respect to the inner product*

$$(y,u) = \int_I p(x)y(x)u(x)dx$$

*Here, in (ii), at an end point at infinity, the difference is to be understood in the sense of limits.*

**Proof** : Let us see how we can recover the weight function $p$ and the boundary conditions in (iii) by demanding that $L$ should be self-adjoint for the inner product

$$(y,u) = \int_I p(x)y(x)u(x)dx$$

for a function $p$ still to be determined .

Let $\alpha, \beta$ be the end points of $I$ . So

$$(Ly,u) = \int_\alpha^\beta p(ay'' + by' + cy)u\, dx \, .$$

Using integration by parts, we have

$$\int_\alpha^\beta (pau)y^{(k)}dx = (pau)y^{(k-1)}\Big|_\alpha^\beta - \int_\alpha^\beta (pau)'y^{(k-1)}dx \, ,$$

so

$$(Ly,u) = (pau)y' - (pau)'y\Big|_\alpha^\beta + \int_\alpha^\beta (pau)''y\,dx + (pbu)y\Big|_\alpha^\beta - \int_\alpha^\beta (pbu)'y\,dx + \int_\alpha^\beta (pcu)y\,dx$$

$$= (pau)y' - (pau)'y + (pbu)y\Big|_\alpha^\beta + \int_\alpha^\beta y((pau)'' - (pbu)' + (pcu))dx$$

This will equal

$$(y,Lu) = \int_\alpha^\beta py(au'' + bu' + cu)\,dx$$

if the boundary terms

$$(pau)y' - (pau)'y + (pbu)y\Big|_\alpha^\beta = 0 \qquad\qquad (*)$$

and

$$\int_\alpha^\beta y((pau)'' - (pbu)' + (pcu))dx = \int_\alpha^\beta yp(au'' + bu' + cu)\,dx$$

The second equality will certainly hold if we restrict $u$ to the class of functions on which the following equation holds:

$$(pau)'' - (pbu)' = pau'' + pbu' \, , \text{ for all } u \, .$$

This requires

$$(pa)''u + 2(pa)'u' + (pa)u'' - (pb' + p'b) = pau'' + pbu' \, .$$

This simplifies to $[(pa)'' - (pb)']u + 2(pa)'u' = 2pbu'$.

Equating coefficients of $u$ and $u'$ on both sides, we get the differential equations for p:

$$(pa)'' - (pb)' = 0 \text{ and } (pa)' = pb \, ,$$

so in fact we need only the equation

$$(pa)' = pb \qquad\qquad (**).$$

The requirement (*) that the boundary terms should vanish now simplifies to

$(pau)y' - (pa)'uy - (pa)u'y + (pbu)y \big|_\alpha^\beta = (pau)y' - (pa)u'y \big|_\alpha^\beta = 0$

or,

$pa(uy' - u'y) = 0$.

The differential equation for the weight is $a'p + ap' = pb$, which integrates to

$$p = e^{\int \frac{(b-a')}{a}dx} = \frac{1}{|a|}e^{\int \frac{b}{a}dx}$$

***Examples***:

1) **Legendre Polynomials**

Consider the eigenvalue problem

$$\left(1 - x^2\right)y'' - 2xy' = \lambda y.$$

Let $L$ be the operator defined by $L(y) = \left(1 - x^2\right)y'' - 2xy'$. The weight function here is $p(x) = \left(\dfrac{1}{\left|1 - x^2\right|}\right)e^{\int \frac{-2x}{1-x^2}dx} = 1$. So in the interval $\left[-1, 1\right]$, $p(x)a(x)$ is $1 - x^2$.

For every non-negative integer $n$, the operator $L$ maps the vector space $P_n$ of polynomials of degree at most $n$ into itself. As $p(x)a(x)$ vanishes at the boundary points, there are no boundary conditions required to make $L$ self adjoint. So $L$ must have a basis of eigenvectors in $P_n$. As $L$ maps every $P_m$ to itself for all $m \le n$ there must be up to a scalar a unique polynomial of every degree $n$ which is an eigenfunction of $L$. We normalize this polynomial so that the coefficient of its leading term is $1$. The corresponding eigenvalue is therefore given by the coefficient of $x^n$ in the expression

$$\left(1 - x^2\right)\left(n(n-1)x^{n-2} + \cdots\right) - 2x\left(nx^{n-1} + \cdots\right).$$

So, it is equal to $-n(n-1) - 2n = -n(n+1)$.

We have therefore recovered the differential equation for the Legendre polynomial: namely $\left(1 - x^2\right)y'' - 2xy' + n(n+1)y = 0$.

## 2) Laguerre Polynomials

Here the eigenvalue problem is

$$xy'' + (1-x)y' = \lambda y.$$

The operator $L$ is given by $L(y) = xy'' + (1-x)y'$. The weight function is

$$p(x) = \left(\frac{1}{|x|}\right)e^{\int \frac{1-x}{x}dx} = e^{-x}.$$ Here the inner product is defined for all functions on

$[0, \infty)$ for the weight $p(x) = e^{-x}$. As $\int\limits_{0}^{\infty} x^n e^{-x} dx = n!$, all polynomials have finite

norm. As $xp(x)$ vanishes at $x = 0$, the operator $L$ will be self adjoint on

the space of all functions with finite norm which are such that $xe^{-x}(uy' - u'y)$

vanishes at $\infty$ for any two functions $u, y$ in this space. As in the first example,

for every non-negative integer $n$, there must be a polynomial of degree $n$ which must be an eigenfunction of $L$. The corresponding eigenvalue is given by the coefficient of $x^n$ in $x(n(n-1)x^{n-2} + \cdots) + (1-x)(nx^{n-1} + \cdots)$.

So, it is equal to $-n$. The differential equation for Laguerre polynomial is

therefore $xy'' + (1-x)y' + ny = 0$

## 3) Hermite Polynomials

The Hermite differential equation is $y'' - 2xy' + \lambda y = 0$. The weight function is

$$p(x) = e^{\int -2x dx} = e^{-x^2}.$$ Thus all polynomials have finite norm relative to this

weight. By the same consideration as in the previous examples there must be a

polynomial of every degree $n$ which is an eigenfunction of $L$ with eigenvalue

$-2n$.

## 4) Confluent hypergeometric equation

This is the equation $xy'' + (c-x)y' - \lambda y = 0$. The operator here is

$L(y) = xy'' + (c-x)y'$. The weight function is therefore $p(x) = |x|^{c-1}e^{-x}$. $L$ maps

polynomials to polynomials and all polynomials have finite norm relative to

$p(x)$ if $c > 1$. As above $L$ maps the space $P_n$ to itself, so by similar considerations as in the previous examples, there must be a polynomial of every degree $n$ which is an eigenfunction of $L$. The corresponding eigenvalue is $-n$. These polynomials are therefore solutions of the differential equation $xy'' + (c-x)y' + ny = 0$.

## 5) Chebyschev Polynomials

These are eigenfunctions of the equations

$$L(y) = \left(1-x^2\right)y'' - xy' \text{ and } L(y) = \left(1-x^2\right)y'' - 3xy'$$

We will discuss only the first case, as the other is similar. The weight function is

$p(x) = \dfrac{1}{\sqrt{(1-x^2)}}$. The singularities at $x = 1, -1$ are not essential singularities.

The reason is that for any continuous function $f$ ,the integral $\displaystyle\int_{-1}^{1} \dfrac{f(x)dx}{\sqrt{(1-x^2}}$ is finite, as one sees by the substitution $x = \cos(\theta)$.

Also, the product of the leading term and the weight is $\sqrt{(1-x^2)}$, so the operator $L$ is self-adjoint on the interval $[-1,1]$ as the contribution from the boundary terms vanishes. Thus, exactly as for the case of Legendre polynomials, there is, upto a constant, exactly one polynomial of degree $n$ which is an eigen function of the operator $L$. The corresponding eigen value is $-n^2$ and these polynomials are solutions of the equation $\left(1-x^2\right)y'' - xy' + n^2 y = 0$

## 6) Jacobi Polynomials

These are eigen functions of the operators $L(y) = (1-x^2)y'' + (\alpha x + \beta)y'$.

Here the weight $p(x)$ can be singular at the end points, but this singularity is inessential, in the sense that $\displaystyle\int_{-1}^{1} p(x)f(x)dx$ is finite for all polynomials $f(x)$, if $\alpha < \beta < -\alpha$. The corresponding differential equation is $(1-x^2)y'' + (ax+b)y' + (n(n-1) - n\alpha)y = 0$: see the Remarks at the end of the following section for details. The Legendre and Chebyschev polynomials are special cases, corresponding to the values $\alpha = -1, -2, -3$ and $\beta = 0$.

## 3. CANONICAL FORMS OF 2$^{\text{ND}}$ ORDER EQUATIONS
## WITH POLYNOMIAL SOLUTIONS

We first summarize Bochner's method, as described in [11]. Bochner defines Sturm-Liouville Polynomial Systems as a sequence of polynomials $P_0(x), ...., P_n(x), ....$, where $P_n(x)$ is of degree $n$ and each of these polynomials is a solution of the differential equation $a(x)y'' + b(x)y' + c(x)y + \lambda y = 0$, for suitable eigenvalues $\lambda = \lambda_0, \lambda_1, .....$

He sets himself the task of determining all these polynomial systems. By considerations of degrees, one sees that the degrees of $a(x), b(x), c(x)$ are at most 2, 1 and 0, respectively. By linear transformation of variables, the normal forms for the possible differential equations are obtained. In each of these cases, the method of undetermined coefficients leads to two-term recursion formulas for the coefficients and one can then decide which eigen values $\lambda = \lambda_k$ give polynomial solutions .

Let us now explain the main difference in the method of this paper with Bochner's method. Although Bochner restricts himself to 2$^{\text{nd}}$ order equations with polynomial coefficients, if one considers a k$^{\text{th}}$ order linear differential equation with rational coefficients, and assumes that there are (k+1) polynomial solutions of different degrees, then one of the eigenvalues $\lambda$ must be non-zero. The reason is that a homogeneous k$^{\text{th}}$ order linear differential equation cannot have more than k linearly independent solutions. Keeping in mind that the degree of a sum of rational functions is at most the degree of each of the summands, we see that the degree of the coefficient term $a_k(x)$ in the term $a_k(x)y^{(k)}$ must be at most $k$ .

Coming back to equations of order 2, if for the operator $L(y) = a(x)y'' + b(x)y' + c(x)y$ there are more than two polynomials of different degrees which are eigenfunctions of $L$, then the degrees of $a(x), b(x), c(x)$ are, at most

2, 1 and 0, respectively. Therefore, $L$ will map the vector space of all polynomials of degree at most $n$ into itself.

From the formulas in Section 2, the operator $L$ would be self adjoint if there is no contribution from the boundary terms: this is ensured if the product $a(x)p(x)$ vanishes at the end points of the interval – finite or infinite- on which the natural weight function $p(x)$ associated to $L$ is finite on the entire interval.

The finiteness of the weight function at the end points then determines the differential equation and finiteness of the norm of polynomials ensures that manipulations as in the proposition of the last section are legitimate. The operator $L$ will then be self-adjoint and it will operates on the vector space of all polynomials of degree at most $n$ for every non-negative integer $n$. As $L$ has a basis of eigenvectors in any finite dimensional subspace on which it operates, we see that there will be, up to a constant, a unique polynomial of degree n, which will be an eigenfunction of $L$, and the corresponding eigenvalues would therefore be determined from the form of the equation. Since these eigenvalues are distinct for different degrees, these polynomials would automatically be orthogonal.

We now determine the operators $L$ from the requirements that

(1): leading term $a(x)$ is non-zero and of degree at most 2, degree of $b(x)$ is 1 at most and $c(x)$ is a constant

(2): the natural weight function associated to $L$ is finite on the interval $I$ where $a(x)$ does not change sign

(3): $a(x)p(x)$ vanishes at the end points of $I$ and, in case there is an end point at infinity, the product $a(x)p(x)$P(x) should vanish at infinity for all polynomials P(x)

(4) All polynomials should have finite norm on $I$ with the weight $p(x)$.

***Case I:*** The polynomial $a(x)$ has two distinct real roots.

By a linear change of variables and scaling we may assume that the roots are $1$ and $-1$. Assuming that $a(x)$ is non-negative in the interval $[-1,1]$, we have $a(x) = 1 - x^2$. Let $b(x) = \alpha x + \beta$ so

$$\frac{b(x)}{a(x)} = \frac{\alpha x + \beta}{(1-x)(1+x)} = \frac{\frac{\beta+\alpha}{2}}{1-x} + \frac{\frac{\beta-\alpha}{2}}{1+x}.$$

So the weight $p(x)$ is

$$p(x) = \frac{1}{1-x^2} e^{\int \left( \frac{\frac{\beta+\alpha}{2}}{1-x} + \frac{\frac{\beta-\alpha}{2}}{1+x} \right) dx} = \frac{(1+x)^{\frac{\beta-\alpha-2}{2}}}{(1-x)^{\frac{\beta+\alpha+2}{2}}}$$

The weight is obviously finite in the interval $(-1,1)$. For $p(x)$ to be finite at the end points, we must have $\beta - \alpha - 2 \geq 0$ and $\beta + \alpha + 2 \leq 0$. Thus $2 + \alpha \leq \beta \leq -2 - \alpha$ and $\alpha \leq -2$. The critical case $\alpha = -2$ gives $\beta = 0$ and the corresponding differential equation is the Legendre differential equation.

***Case II:*** The polynomial $a(x)$ has repeated real roots.

In this case we can assume that $a(x) = x^2$. Let $b(x) = \alpha x + \beta$. The weight function is now $p(x) = \frac{1}{x^2} e^{\int \frac{\alpha x + \beta}{x^2} dx} = \frac{|x|^\alpha}{x^2} e^{-\beta/x}$.

If $\beta > 0$, then $\lim_{x \to 0^+} \frac{|x|^\alpha}{x^2} e^{-\beta/x}$ is finite and $\lim_{x \to 0^-} \frac{|x|^\alpha}{x^2} e^{-\beta/x}$ is infinite. If $\beta < 0$, then $\lim_{x \to 0^+} \frac{|x|^\alpha}{x^2} e^{-\beta/x}$ is infinite and $\lim_{x \to 0^-} \frac{|x|^\alpha}{x^2} e^{-\beta/x}$ is finite. If $\beta = 0$, then the weight function is $|x|^{\alpha-2}$ which is finite at zero only if $\alpha \geq 2$ but in this case the integral from $-\infty$ to $0$ or from $0$ to $\infty$ of $|x|^{\alpha-2}$ is infinite. So for a finite weight, we must take either $\beta < 0$ and I to be $[0,\infty)$ or $\beta > 0$ and the interval $I$ to be $(-\infty,0]$. Since the discussion for both these cases is identical, we limit ourselves only to the case $\beta > 0$ and interval $I = [0,\infty)$. The weight function is now $p(x) = \frac{x^\alpha}{x^2} e^{-\beta/x}$.

Since $\int_0^\infty e^{-x} x^p dx$ is finite only if $p > -1$ we see that the integral $\int_0^\infty x^k e^{-\left(\frac{\beta}{x}\right)} dx$ is finite only if $k < -1$. So polynomials do not have a finite norm for the weight $p(x) = \dfrac{x^\alpha}{x^2} e^{-\beta/x}$ and this case therefore vacuous.

***Case III:*** The polynomial $a(x)$ is linear.

In this case we can take $a(x) = x$. Let $b(x) = \alpha x + \beta$. In this case the weight function is

$$p(x) = \frac{1}{|x|} e^{\int \frac{\alpha x + \beta}{x} dx} = |x|^{\beta - 1} e^{\alpha x}.$$ This is finite at zero if and only if $\beta \geq 1$. Since $\int_0^\infty e^{\alpha x} x^\varepsilon dx$,

where $\varepsilon > 0$, is finite only if $\alpha \leq 0$ we see that we can not take the interval $I$ from $-\infty$

to $\infty$. Without loss of generality we can take this to be the interval $[0, \infty)$. So the weight function is now $p(x) = x^{\beta - 1} e^{\alpha x}$ with $\alpha < 0$ and $\beta \geq 1$. All polynomials have finite norm with respect to this weight and for all polynomials $P(x)$ the product $P(x) p(x)$ vanishes at $0$ and $\infty$. Therefore the equation $xy'' + (\alpha x + \beta) y' + \lambda y = 0$ has polynomial solutions for every degree $n$. The corresponding eigenvalue is $\lambda = -\alpha n$.

***Case IV:*** $a(x) = 1$

In this case $L(y) = y'' + (\alpha x + \beta) y' + \gamma y$. The weight is $p(x) = e^{\frac{\alpha x^2}{2}} e^{\beta}$. So $\alpha$ must be negative, for the product $P(x) p(x)$ to vanish at the end points of the interval $I$ for all polynomials $P(x)$, and therefore $I$ must be $(-\infty, \infty)$.

**Remarks**

- The case of $a(x)$ with no real roots does not arise, because of the requirements
  (3) and (4) above which a weight function must satisfy.

- If one allows non-essential singularities of the weight function- as in the
  example of Chebyschev polynomials- example 5, Section 2- one will recover a
  wider class of equations with polynomial solutions- classes with minimal degree
  of admissible polynomials $\geq 1$ indeed .

We give details here for one important case, namely, the Jacobi differential equation. First note that for any differentiable function $f$ with $f'$ continuous, the integral $\int_0^e \dfrac{f(x)}{x^\alpha}dx$ is finite if $\alpha < 1$ - as one sees by using integration by parts.

Consider the equation $(1-x^2)y'' + (ax+b)y' + \lambda(y) = 0$. As in *Case I* above, the weight function $p(x)$ for the operator $L(y) = (1-x^2)y'' + (ax+b)y'$ is $p(x) = \dfrac{1}{1-x^2}e^{\int\left(\frac{\frac{\beta+\alpha}{2}}{1-x} + \frac{\frac{\beta-\alpha}{2}}{1+x}\right)dx} = \dfrac{(1+x)^{\frac{\beta-\alpha-2}{2}}}{(1-x)^{\frac{\beta+\alpha+2}{2}}} = \dfrac{1}{(1-x)^{\frac{\beta+\alpha+2}{2}}(1+x)^{\frac{-\beta+\alpha+2}{2}}}$ .

So $\int_{-1}^{1} p(x)f(x)dx$ would be finite if $\beta + \alpha < 0$ and $-\beta + \alpha < 0$, that is, if $\alpha < \beta < -\alpha$ . Moreover, $(1-x^2)p(x) = (1-x)^{\frac{-(\beta+\alpha)}{2}}(1+x)^{\frac{\beta-\alpha}{2}}$ vanishes at the end-points -1 and 1. Therefore, $L$ is a self-adjoint operator on all polynomials of degree $n$ and so, there must be, up to a scalar, a unique polynomial which is an eigen function of $L$ for eigen value $-n(n-1) + n\alpha$ .

So these polynomials satisfy the equation $(1-x^2)y'' + (ax+b)y' + (n(n-1) - n\alpha)y = 0$, for any negative real number $\alpha$ . This is an alternative form of the Jacobi and Gauss differential equations.

## 4. DETERMINING EQUATIONS FOR HIGHER ORDER EQUATIONS

In this section, we give the determining equations for the weight and constraints on the coefficients which make equation of order 4 and 3 into a self adjoint or an anti-self adjoint system respectively, leaving the details for this and other higher order equations to [1]. If one wants self adjoint equations in all dimensions, one has to consider differential equations with complex coefficients and Hermitian inner products.

The key to all this is the case of a linear ODE of order one and the behavior of the weight function at singular points of the equation. A detailed analysis of this follows. We conclude this section by giving explicit examples of third and fourth order equations.

***Operators of order 1***

Consider the operator $L(y) = a(x)y' + b(x)y$. $L$ is anti self adjoint if $(Ly, u) = -(y, Lu)$, the inner product is taken for a weight function $p(x)$ yet to be determined. As in the two-dimensional case, the differential equation for $p(x)$ is $(ap)' = 2bp$ which integrates to $p(x) = \dfrac{1}{|a|}e^{\int \frac{2b}{a}dx}$ and the boundary conditions needed for $L$ are $a(x)p(x)u(x)y(x)\big|_I = 0$.

Clearly if $a(x)$ has no zero in an interval $I$ then $L$ will be anti self adjoint on any linear space $C$ of functions $C$ for which the following hold:

- $C$ is invariant under $L$

- for every function $y$ in $C$, the integral of $py^2$ over $I$ is finite

- for any functions $u$ and $y$ in $C$ the difference of $a(x)p(x)u(x)y(x) = apuy$ evaluated at the end points of $I$ is zero.

So the interesting case is the singular case where $a(x)$ has a zero in $I$. Assume $a(x) = 0$ at $x = 0$. Suppose $a(x) = x^\alpha g(x)$ with $\alpha$ a positive integer, $g(0) \neq 0$ and $b(x) = x^\beta f(x)$ with $\beta \geq 0$ and $f(0) \neq 0$. So $\dfrac{b(x)}{a(x)} = x^{\beta-\alpha}h(x)$ with $h(0) \neq 0$. Therefore near zero, $p(x) = \dfrac{1}{|x^\alpha||g(x)|}e^{\int 2x^{\beta-\alpha}h(x)dx}$. So near $x = 0$, $p(x)$ is finite if and only if $\dfrac{1}{|x^\alpha|}e^{\int 2x^{\beta-\alpha}h(x)dx}$ is finite near $x = 0$. Now $\lim\limits_{x\to 0} h(x) = h(0) \neq 0$, we have the following possibilities.

(i). If $h(0) > 0$ then for $x$ sufficiently near zero we $\dfrac{1}{2}h(0) < h(x) < \dfrac{3}{2}h(0)$.

(ii). If $h(0) < 0$ then for $x$ sufficiently near zero we $\dfrac{3}{2}h(0) < h(x) < \dfrac{1}{2}h(0)$

In case (i) we have the estimate $\dfrac{1}{|x^\alpha|}e^{\int x^{\beta-\alpha}h(0)dx} < \dfrac{1}{|x^\alpha|}e^{\int 2x^{\beta-\alpha}h(x)dx} < \dfrac{1}{|x^\alpha|}e^{\int 3x^{\beta-\alpha}h(0)dx}$.

In case (ii) we have the estimate $\dfrac{1}{|x^\alpha|}e^{\int 3x^{\beta-\alpha}h(0)dx} < \dfrac{1}{|x^\alpha|}e^{\int 2x^{\beta-\alpha}h(x)dx} < \dfrac{1}{|x^\alpha|}e^{\int x^{\beta-\alpha}h(0)dx}$.

Using the facts that for the limit $\displaystyle\lim_{x \to 0} \frac{e^{(\lambda x^{\beta})}}{|x^{\alpha}|}$

- the right hand limit exists and equals zero only if $\lambda < 0$ and $\beta < 0$ and the left hand limit exists and equals zero only if $(-1)^{\beta}\lambda < 0$ and $\beta < 0$

- the two sided limit exists and equals zero only if $\lambda < 0$ and $\beta$ is even and negative

one can determine the allowable weights and the type of the interval (closed or half closed or $(-\infty, \infty)$).

### *Operators of order 3*

Consider the operator $L(y) = a_3(x)y''' + a_2(x)y'' + a_1(x)y' + a_0(x)y$. The condition that $L$ is anti self adjoint gives the following determining equations for the weight $p(x)$, the coefficient functions and the required boundary conditions:

$$\left(a_3 p\right)' = \frac{2}{3}a_2 p$$

$$\left(a_2 p\right)'' - 3\left(a_1 p\right)' + 6\left(a_0 p\right) = 0$$

At the boundary, for any functions $u, y$ the difference of the following expression evaluated at the end points of $I$ must equal zero:

$$a_3 p\left[uy'' - u'y' + u''y\right] + \frac{a_2 p}{3}\left[uy' + u'y\right] - \frac{1}{3}\left((a_2 p)' - 3a_1 p\right)uv$$

### *Operators of order 4*

Consider the operator $L(y) = a_4(x)y^{(iv)} + a_3(x)y''' + a_2(x)y'' + a_1(x)y' + a_0(x)y$. The condition that $L$ is self adjoint gives the following determining equations for the weight $p(x)$, the coefficient functions and the required boundary conditions:

$$\left(a_4 p\right)' = \frac{1}{2}a_3 p \qquad\qquad\qquad (4.1)$$

$$\left(a_3 p\right)'' - 2\left(a_2 p\right)' + 2\left(a_1 p\right) = 0 \qquad\qquad (4.2)$$

At the boundary, for any functions $u, y$ the difference of the following expression evaluated at the end points of $I$ must equal zero:

$$a_4 p \left[ uy''' - u'y'' + u''y' - u'''y \right] + \frac{a_3 p}{2} \left[ uy'' - u''y \right] - \frac{1}{2} \left( (a_3 p)' - 2a_2 p \right) \left[ uy' - u'y \right] \quad (4.3)$$

We conclude this section by deriving an example of a fourth order linear differential operator which has polynomial eigenfunctions of every degree $n$. A special case of this example is given in [5], which also gives other examples corresponding to families of classical orthogonal polynomials. These examples can be also be obtained by choosing suitable weights and leading terms as in the following example. We will give a general classification in a subsequent paper.

<u>Example</u>

If we take $p = 1$ and $a_4(x) = \left( 1 - x^2 \right)^2$ then using Equations (4.1) and (4.2) we have

$a_3(x) = -8x(1 - x^2)$ and $a_2'(x) = 24x + a_1(x)$ and $a_0(x)$ arbitray.

The expression (4.3) becomes

$$\left( 1 - x^2 \right)^2 \left[ uy''' - u'y'' + u''y' - u'''y \right] - 4x \left( 1 - x^2 \right) \left[ uy'' - u''y \right] + \left( a_2 + 4 - 12x^2 \right) \left[ uy' - u'y \right] \quad (4.4)$$

Choosing $a_2(x)$ to be an even function we see that the difference at the end points of the interval $I = \left[ -1, 1 \right]$ of the expression in (4.4) is zero.

In particular, for $a_2(x) = 8$, the equation

$$\left( 1 - x^2 \right)^2 y^{(iv)} - 8x \left( 1 - x^2 \right) y''' + 8y'' - 24xy' - \lambda y = 0$$

must have polynomial solutions of degree $n$ for the eigenvalue $\lambda = n \left( (n-1)(n-2)(n+5) - 24 \right)$.